  \documentclass[preprint,12pt]{article}
\usepackage{amssymb}
 \usepackage{amsthm}
\usepackage{latexsym}
\usepackage{amsmath}
\usepackage{amssymb}
\usepackage{color}

\newcommand{\udots}{\mathinner{\mskip1mu\raise1pt\vbox{\kern7pt\hbox{.}}
\mskip2mu\raise4pt\hbox{.}\mskip2mu\raise7pt\hbox{.}\mskip1mu}}
\usepackage{amsmath}

\newtheorem{theorem}{Theorem}[section]

\newtheorem{lemma}{Lemma}[section]
\newtheorem{remark}{Remark}[section]

\begin{document}
\title{\bf \Large { The Bessel function expression of characteristic function}}
{\color{red}{\author{
\normalsize{Chuancun Yin \;\;Hua Dong}
\thanks{Corresponding author.}\\
{\normalsize\it  School of Statistics and Data Science,  Qufu Normal University}\\
\noindent{\normalsize\it Qufu 273165, Shandong,  China}\\
e-mail:   sddh1978@126.com}}}
\maketitle
\vskip0.01cm
\centerline{\large {\bf Abstract}}
 In this paper, we   give a unified method to derive the classical characteristic functions of all elliptical  and  related
distributions in terms of  Bessel functions.  The approach is based on the stochastic representation of elliptical random variable
and the characteristic  function of
uniform distribution on the unit sphere surface in $\Bbb{R}^n$.  In particular, we present the  simple closed form of characteristic functions for  commonly used distributions such as  multivariate $t$,  Pearson Type II, Pearson Type VII, Kotz type  and  Bessel distributions. Some extensions are also investigated.

\medskip

\numberwithin{equation}{section}
\section{Introduction}\label{intro}

The characteristic function (CF) of a random vector ${\bf X}$ is the function $\psi_{\bf X}: {\Bbb R}^n \rightarrow {\Bbb C}$ defined as
$\psi_{\bf X}({\bf t})=E(e^{i{\bf t}'{\bf X}}),\; {\bf t}\in\Bbb{R}^n.$
Characteristic functions  play an important role in probability and statistics. In recent years, there has been a growing interest in the  multivariate  elliptical distributions and related  multivariate   distributions such that the skew-elliptical distributions, location-scale mixture of elliptical  distributions, and so on.
The characteristic functions of many    multivariate   distributions have been derived by many authors. For instance,  a general
treatment for  multivariate  elliptical distributions  has been given in  Fang, Kotz and Ng (1990).  Arellano-Valle and Genton (2005) studied the moment generating function for fundamental skew-symmetric distributions,  Arellano-Valle and Azzalini (2006)  studied the moment generating function for unified skew-normal distributions. Kim and Genton (2011) concentrated on the characteristic functions of scale mixtures of multivariate skew-normal distributions which include some well-known distributions, for example, the skew-normal and skew-$t$.   Vilca, Balakrishnan and  Zeller (2014) investigated the multivariate   skew-normal generalized hyperbolic distribution and derived its moment generating functions.  Some recent results show that the characteristic functions are very useful in the study of stochastic orderings and  deriving the moments  for multivariate skew-elliptical distributions, see, e.g., Yin (2021),  Jamali et al. (2020),  Pu et al. (2021), and Yin and  Balakrishnan (2024), and so on.

In the  existing literature,  most of  characteristic functions   are derived
based on  contour integrals which unfamiliar to statistics students, or by solving  ordinary differential equations,  some  expressions  are too complex so as to limit the application in practice, and even some expressions are wrong. The purpose of this short note is
to complement the existing literature with a new simple and direct derivation  for the
 multivariate elliptical distributions and their location-scale mixtures, which include some well-known
distributions  such as the multivariate normal, multivariate $t$,   multivariate  Kotz-type,  multivariate Pearson type VII, multivariate stable law, Cauchy, logistic, Laplace and so on. The advantage of this method is that the calculation process is simple, the contour  integrals are avoided, and the expressions are  more concise.

The paper is organized as follows.  In Section 2, we recall some important concepts of  elliptical symmetric distributions, several special functions  including  generalized hypergeometric series and  Bessel functions and characteristic  function of the  uniform distribution on the unit sphere surface in $\Bbb{R}^n$. Section 3 is the main results, in which we present a unified   expression  for the characteristic functions  of  multivariate elliptical distributions     in terms of Bessel functions, and provides simple closed form of characteristic functions for  commonly used distributions. Some extensions are given in Section 4.
Finally, Section 5 provides the conclusions.

\section{ Preliminaries}
In this section, we first present a brief overview of  elliptical  distributions,  some special functions  such as the generalized hypergeometric series and the  Bessel functions, and an  alternative expression for  characteristic  function of the uniform distribution on the unit sphere surface in $\Bbb{R}^n$.
\subsection{Elliptical distributions}
The class of  elliptical distributions provides a generalization of the multivariate normal distributions. This class includes various distributions such as symmetric Kotz type distribution, symmetric multivariate Pearson type VII and symmetric multivariate stable law, multivariate Student-$t$, Cauchy, logistic, Laplace and so on.
They are useful to modeling multivariate random phenomena which have heavier tails than the normal as well as having some skewness. Such a rich
class of distributions can be used to model multivariate regression
problems with skew-elliptical error structure.

An $n \times 1$ random vector ${\bf X}= (X_1, X_2,\cdots, X_n)'$ is said to have an elliptical
 distribution if its characteristic function is $e^{i{\bf t}'{\boldsymbol \mu}}\phi({\bf t}'{\bf \Sigma}{\bf t})$ for all ${\bf t}\in \Bbb{R}^n $,  where  $\phi$ is called the characteristic generator satisfying $\phi(0)=1$,
$\boldsymbol{\mu}$ ($n$-dimensional vector) is its location parameter    and  $\bf{\Sigma}$ ($n\times n$ positive semi-definite matrix) is its dispersion matrix.  We shall write
${\bf{X}}\sim ELL_n ({\boldsymbol \mu},{\bf \Sigma},\phi)$. It is well known that
$\bf X$ admits the stochastic representation
\begin{equation}
{\bf X}={\boldsymbol \mu}+R{\bf A}'{\bf U}^{(n)}, \label{(2.1)}
\end{equation}
where ${\bf A}$  is a square matrix such that ${\bf A}'{\bf A}= {\bf \Sigma}$, ${\bf U}^{(n)}$ is uniformly distributed on the unit sphere surface in $\Bbb{R}^n$,  $R\ge 0$   is the random variable with $R \sim F$ in $[0, \infty)$ called the generating variate and $F$ is called the generating distribution function, $R$ and  ${\bf U}^{(n)}$ are  independent.
In general,   an elliptically distributed random vector ${\bf{X}}\sim ELL_n ({\boldsymbol \mu},{\bf \Sigma},\phi)$ does not necessarily
possess a density. However, if density of ${\bf X}$ exists it must be of the form
\begin{equation}
f({\bf x})=c_n|{\bf \Sigma}|^{-\frac{1}{2}}g(({\bf x}- {\boldsymbol \mu})^{T}{\bf \Sigma}^{-1}({\bf x-{\boldsymbol \mu} })), \;{\bf x}\in \Bbb{R}^n,
\end{equation}
for some non-negative function $g$ satisfying the condition
$$\int_0^{\infty}z^{\frac{n}{2}-1} g(z)dz<\infty, $$
and a normalizing constant $c_n$ given by
\begin{equation}
c_n=\frac{\Gamma(\frac{n}{2})}{\pi^{\frac{n}{2}}}\left(\int_0^{\infty}z^{\frac{n}{2}-1} g(z)dz\right)^{-1}.
\end{equation}
  The function $g$ is called the density generator. One sometimes
writes ${\bf X} \sim ELL_n ({\boldsymbol \mu},{\bf \Sigma},g)$ for the $n$-dimensional elliptical distributions generated from the
function $g$. In this case $R$ in (2.1) has the pdf given by
\begin{equation}
h_R(v)=c_n \frac{2\pi^{\frac{n}{2}}}{\Gamma(\frac{n}{2})}v^{n-1}g(v^2), v\ge 0.
\end{equation}
 A comprehensive review of the properties and characterizations of elliptical distributions can be found  in  Cambanis et al. (1981),  Fang, Kotz and Ng (1990); see also some recent papers, Zuo, Yin and  Balakrishnan (2021), Wang and Yin (2021) and   Yin, Wang and Sha (2022).

\subsection{Special functions}
The following special functions and mathematical results will be useful in our
analyses. The details  can be found in  Gradshteyn and  Ryzhik (2007) and Slater (1960). The series
$${}_1F_2(\alpha_1; \beta_1,\beta_2;z)=\sum_{k=0}^{\infty}\frac{(\alpha_1)_k}{(\beta_1)_k(\beta_2)_k}\frac{z^k}{k!},$$
is called a generalized hypergeometric series of order $(1, 2)$, where $(\alpha)_k$ and
$(\beta)_k$ represent Pochhammer symbols and
$(x)_k = x(x+1)\cdots(x +k-1)$, and $(x)_0 = 1$.
The series given by
$${}_2F_1(\alpha,\beta; \gamma;z)\equiv F(\alpha,\beta; \gamma;z)=
\sum_{k=0}^{\infty}\frac{(\alpha)_k(\beta)_k}{(\gamma)_k}\frac{z^k}{k!},$$
is called a generalized hypergeometric series of order $(2, 1)$. The series
$${}_1F_1(\alpha; \gamma;z)\equiv F(\alpha, \gamma, z)=
\sum_{k=0}^{\infty}\frac{(\alpha)_k}{(\gamma)_k}\frac{z^k}{k!},$$
is called a  confluent hypergeometric function or Kummer's function. The series
$${}_0F_1(\gamma;z)=\sum_{k=0}^{\infty}\frac{1}{(\gamma)_k}\frac{z^k}{k!},$$
is called a a generalized hypergeometric series of order $(0, 1)$.
Bessel function of the first kind, $J_{\nu}(x)$, is defined as
$$J_{\nu}(x)=\sum_{k=0}^{\infty}(-1)^k\frac{1}{k!\Gamma(\nu+k+1)}\left(\frac{x}{2}\right)^{\nu+2k},$$
or
$$J_{\nu}(x)=\frac{x^{\nu}}{2^{\nu}\Gamma(\nu+1)} {}_0F_1\left(\nu+1;-\frac{x^2}{4}\right).$$
Modified Bessel function of the first kind, $I_{\nu}(x)$,  is defined as
$$I_{\nu}(x)=\sum_{k=0}^{\infty}\frac{1}{k!\Gamma(\nu+k+1)}\left(\frac{x}{2}\right)^{\nu+2k}.$$
Modified Bessel function of the  second kind (also called the MacDonald function) $K_{\nu}(x)$     with order $\nu$ may be defined
by the following integral:
$$K_{\nu}(x)=\left(\frac{2}{x}\right)^{\nu}\frac{\Gamma(\nu+\frac12)}{\sqrt{\pi}}\int_0^{\infty}\frac{\cos (xu)}{(1+u^2)^{\nu+\frac12}}du, x>0,\nu>-\frac12.$$
 $I$ and $K$  have asymptotic properties: For $\nu>0$, we have
$$x^{-\nu}I_{\nu}(x)\to \frac{1}{2^{\nu}\Gamma(1+\nu)},\;x^{\nu}K_{\nu}(x)\to 2^{\nu-1}\Gamma(\nu), {\rm as}\; x\to 0.$$

\subsection{Characteristic  function of ${\bf U}^{(n)}$ }
Let ${\bf U}^{(n)}$  be uniformly distributed on the unit sphere surface in $\Bbb{R}^n$
which has been introduced in Section 1.  Let $\Omega_n(||{\bf t}||^2), {\bf t}\in \Bbb{R}^n$ be the characteristic function of ${\bf U}^{(n)}$, where $||{\bf t}||^2={\bf t}'{\bf t}$. The well-known three equivalent forms of $\Omega_n(||{\bf t}||^2)$ can be found in  Fang et al. (1990, p.70). Another result for  characteristic  function of ${\bf U}^{(n)}$
 is due to Schoenberg  (1938) which will  be useful for the rest of the paper. Here we present a simple proof by using an integral formula for $J_{\nu}$.

\begin{lemma} (Schoenberg (1938)).  The characteristic  function of ${\bf U}^{(n)}$ can be expressed as
\begin{equation}
\Omega_n(||{\bf t}||^2)=\Gamma\left(\frac{n}{2}\right)\left(\frac{2}{||{\bf t}||}\right)^{\frac{n-2}{2}}J_{\frac{n-2}{2}}(||{\bf t}||),\; {\bf t}\in \Bbb{R}^n,
\end{equation}
where $J_{\nu}$ is the Bessel function of the first kind of order $\nu$.
\end{lemma}
{\bf Proof} \; It is follows from the proof of  Theorem 3.1 in Fang et al. (1990) that
$$\Omega_n(||{\bf t}||^2)=\frac{\Gamma\left(\frac{n}{2}\right)}{\Gamma\left(\frac{n-1}{2}\right)\sqrt{\pi}}\int_{-1}^{1} e^{i||{\bf t}||u}(1-u^2)^{\frac{n-1}{2}-1}du.$$
By using the formula
$$\int_{-1}^{1} e^{i\mu x}(1-x^2)^{\nu-1}dx=\sqrt{\pi}\Gamma(\nu)\left(\frac{2}{\mu}\right)^{\nu-\frac{1}{2}}J_{\nu-\frac{1}{2}}(\mu),$$
we arrive at the result.

\begin{remark} Using the relationship of   $J_{\nu}$ and ${}_0F_1$ leads the following well known fact (see, e.g.  Fang et al. (1990), P.69):
$$\Omega_n(||{\bf t}||^2)={}_0F_1\left(\frac{n}{2};-\frac{||{\bf t}||^2}{4}\right), {\bf t}\in \Bbb{R}^n.$$
\end{remark}

\numberwithin{equation}{section}
\section { Main results}

 In this section, we derive the CFs of  elliptical  distributions. Moreover, we give some important special cases.
The following  result provided a   concrete  form of the     characteristic generator  function in terms of  the Bessel functions for any elliptical  distribution with density generator $g$.
\begin{theorem}
Suppose  ${\bf X} \sim ELL_n ({\boldsymbol \mu},{\bf \Sigma},g)$. Then, the characteristic function of ${\bf X}$  is given by
\begin{equation}
\psi_{\bf X}({\bf t})=e^{i{\bf t}'{\boldsymbol \mu}}\phi({\bf t}'{\bf \Sigma}{\bf t}),\; {\bf t}\in \Bbb{R}^n
\end{equation}
for some real valued function $\phi$ defined as
\begin{equation}
\phi(u^2)=c_n(2\pi)^{\frac{n}{2}} u^{-\frac{n-2}{2}}\int_0^{\infty}r^{\frac{n}{2}}J_{\frac{n-2}{2}}(ru)g(r^2)dr, \; u\ge 0,
\end{equation}
where $c_n$ is the normalizing constant given by (2.3).
\end{theorem}
{\bf Proof} \; Without loss of generality we assume that ${\boldsymbol \mu}={\bf 0}, {\bf \Sigma}={\bf I_n}$.  For ${\bf t}\in \Bbb{R}^n$, by using the stochastic representation (2.1) and (2.4), we obtain
\begin{eqnarray*}
\phi({\bf t}'{\bf t})&=&E(e^{iR{\bf t}'{\bf U}^{(n)}})=\int_0^{\infty}E(e^{ir{\bf t}'{\bf U}^{(n)}})P(R\in dr)\\
&=&\int_0^{\infty}\Omega_n(r^2||{\bf t}||^2)h_R(r)dr\\
&=&\int_0^{\infty}\Omega_n(r^2||{\bf t}||^2)c_n \frac{2\pi^{\frac{n}{2}}}{\Gamma(\frac{n}{2})}r^{n-1}g(r^2)dr.
\end{eqnarray*}
Therefore, utilizing (2.5) we have
$$\phi({\bf t}'{\bf t})=c_n(2\pi)^{\frac{n}{2}} (||{\bf t}||)^{-\frac{n-2}{2}}\int_0^{\infty}r^{\frac{n}{2}}J_{\frac{n-2}{2}}(r ||{\bf t}||)g(r^2)dr,$$
and we arrive at (3.2). This completes the proof.\\

  The following theorem discusses the inverse transformation problem which was put forward by one of  reviewers.
\begin{theorem}
Suppose  ${\bf X} \sim ELL_n ({\boldsymbol \mu},{\bf \Sigma}, \phi)$ has  a density. Then, the    density generator of ${\bf X}$  is given by
\begin{equation}
 g(r)=\frac{r^{-\frac{n-2}{4}}}{c_n(2\pi)^{\frac{n}{2}}}\int_0^{\infty}u^{\frac{n}{2}}J_{\frac{n-2}{2}}(\sqrt{r}u)\phi(u^2)du, \; r\ge 0,
\end{equation}
where $c_n$ is the normalizing constant given by (2.3).
\end{theorem}
{\bf Proof} \; Using the Hankel transform pair (see, Davies (2002), P.227), the integral transform  (3.2)  has   the reciprocal formula
$$r^{\frac{n}{2}-1} g(r^2)=\frac{1}{c_n(2\pi)^{\frac{n}{2}}}\int_0^{\infty}u^{\frac{n}{2}}J_{\frac{n-2}{2}}(ru)\phi(u^2)du, \; r\ge 0, $$
from which we get (3.3) immediately,  ending the proof.

In the following, we shall using Theorem 3.1 to derive the CFs of  uniform distribution in the unit sphere, multinormal,  multivariate $t$, Pearson Type II, Pearson Type VII, Kotz type and multivariate Bessel.  Of course, most of results are not new, the proposed method offers simplicity and ease of use, making it potentially valuable from a pedagogical perspective.

\subsection{Uniform distribution in the unit sphere}

Let ${\bf V^{(n)}}$ be distributed according to the uniform
distribution in the unit sphere in   $\Bbb{R}^n$ (see Fang et al. (1990), P.74).  ${\bf V^{(n)}}$ has the  stochastic representation
 ${\bf V^{(n)}}=R {\bf U}^{(n)}$. The density of $R$ is given by
 $$f(r)=nr^{n-1}{\bf 1}_{(0\le r\le 1)},$$
 and the density generator of  ${\bf V^{(n)}}$   is given by
$$g(r)=1,\; 0\le r\le 1.$$
Then $$c_n=\frac{\Gamma(\frac{n}{2})}{\pi^{\frac{n}{2}}}\left(\int_0^{1}z^{\frac{n}{2}-1}dz\right)^{-1}=\frac{n\Gamma(\frac{n}{2})}{2\pi^{\frac{n}{2}}}.$$
By using (3.1) and (3.2) we find the characteristic function of ${\bf V^{(n)}}$ given by
 \begin{eqnarray}
\psi_{\bf V^{(n)}}({\bf t})&=&2^{\frac{n-2}{2}}n\Gamma\left(\frac{n}{2}\right)||{\bf t}||^{-\frac{n-2}{2}}\int_0^{1}r^{\frac{n}{2}}J_{\frac{n-2}{2}}(r ||{\bf t}||)dr\nonumber\\
&=&2^{\frac{n-2}{2}}n\Gamma\left(\frac{n}{2}\right)||{\bf t}||^{-\frac{n}{2}}J_{\frac{n}{2}}(||{\bf t}||), \;  {\bf t}\in \Bbb{R}^n,
\end{eqnarray}
where we have used the integral (Gradshteyn and Ryzhik (2007, p.676)):
$$\int_0^1 x^{\nu+1}J_{\nu}(ax)dx=\frac{1}{a}J_{\nu+1}(a).$$
Equivalently,
\begin{eqnarray}
\psi_{\bf V^{(n)}}({\bf t})={}_0F_1\left(\frac{n}{2}+1;-\frac{||{\bf t}||^2}{4}\right), \;  {\bf t}\in \Bbb{R}^n.
\end{eqnarray}
\begin{remark} The characteristic function of ${\bf V^{(n)}}$ is  given in Fang et al. (1990), P.75) by
$$\phi({\bf t}'{\bf t})=\frac{2}{B(\frac12, \frac{n+1}{2})}\int_0^{\infty}\cos({\bf t}'{\bf t} x)(1-x^2)^{\frac{n-1}{2}}dx.$$
 We remark that the $\phi$ above does not  satisfies $\phi(0)=1$. Actually,
\begin{eqnarray}
\phi({\bf t}'{\bf t})=\frac{2}{B(\frac12, \frac{n+1}{2})}\int_0^{1}\cos(||{\bf t}||x)(1-x^2)^{\frac{n-1}{2}}dx.
\end{eqnarray}
By using  formula 8 in Gradshteyn and Ryzhik (2007, p.442)), which states that, for $a>0, u>0, \Re(u)>-\frac12$,
$$\int_0^u (u^2-x^2)^{\nu-\frac12}\cos(ax)dx=\frac{\sqrt{\pi}}{2}\left(\frac{2u}{a}\right)^{\nu}\Gamma\left(\nu+\frac12\right)J_{\nu}(au).
$$
we can prove that (3.3) and (3.5) are equivalent.
\end{remark}

\subsection{Multinormal }
The random vector ${\bf X} = (X_1, X_2,\cdots, X_n)'$ is said to have a multivariate normal distribution  $N_n ({\boldsymbol \mu},{\bf \Sigma})$
if its  density generator is of the form $g_n(u)=\exp{(-\frac{u}{2})}$ and the normalizing constant is given by
$c_n=(2\pi)^{-\frac{n}{2}}.$
Applying (3.2) to $c_n$ and $g_n$ we get
\begin{equation}
\phi(u^2)= u^{-\frac{n-2}{2}}\int_0^{\infty}r^{\frac{n}{2}}J_{\frac{n-2}{2}}(ru) \exp\left(-\frac{u^2}{2}\right)dr.
\end{equation}
By the following  formula   (Gradshteyn and Ryzhik (2007, p.706))
$$\int_0^{\infty} x^{\mu}e^{-ax^2}J_{\nu}(\beta x)dx=\frac{\beta^{\nu}\Gamma(\frac{1}{2}\mu+\frac{1}{2}\nu+\frac{1}{2})}{2^{\nu+1}\alpha^{\frac{1}{2}\mu+\frac{1}{2}\nu+\frac{1}{2}}\Gamma(\nu+1)} {}_1F_1\left(\frac{1}{2}\mu+\frac{1}{2}\nu+\frac{1}{2}; \nu+1, -\frac{\beta^2}{4\alpha}\right),$$
 one  gets that
$$\phi(u^2)={}_1F_1\left(\frac{n}{2}; \frac{n}{2}, -\frac{u^2}{2}\right)=e^{-\frac{1}{2}u^2}.$$
Thus, we get the following  well known  characteristic function of $N_n ({\boldsymbol \mu},{\bf \Sigma})$:
$$\psi_{\bf X}({\bf t})=e^{i{\bf t}'{\boldsymbol \mu}}\exp\left(-\frac12 {\bf t}'{\bf \Sigma}{\bf t}\right),\; {\bf t}\in \Bbb{R}^n.$$

\subsection{ Multivariate $t$ }

The $n$-dimensional random vector ${\bf X}$ is said to have a multivariate generalized $t$ distribution  if its probability density function  is given by (see Fang et al. (1990))
$$f({\bf x})=\frac{c_n}{\sqrt{|{\bf \Sigma}|}}\left(1+ \frac{({\bf x}- {\boldsymbol \mu})^{T}{\bf \Sigma}^{-1}({\bf x-{\boldsymbol \mu}) }}{s}\right)^{-\frac{n+m}{2}},\;  {\bf x}\in \Bbb{R}^n,  s>0, m\in \Bbb{N}_+.$$
If $m =s$, then ${\bf X}$  follows the multivariate $t$ distribution.
Its characteristic function   is given by
\begin{equation}
\psi_{\bf X}({\bf t})=e^{i{\bf t}'{\boldsymbol \mu}}\phi({\bf t}'{\bf \Sigma}{\bf t}),\; {\bf t}\in \Bbb{R}^n,
\end{equation}
where
\begin{equation}
\phi(u^2)=\frac{2^{\frac{n}{2}}\Gamma(\frac{n}{2})}{\int_0^{\infty}z^{\frac{n}{2}-1} g(z)dz}u^{-\frac{n-2}{2}}
\int_0^{\infty}r^{\frac{n}{2}}J_{\frac{n-2}{2}}(ru)g(r^2)dr.
\end{equation}
Here,
$$g(u)=\left(1+\frac{u}{s}\right)^{-\frac{n+m}{2}}.$$
It follows the formula (see  Gradshteyn and Ryzhik (2007, p.678))
\begin{equation}
\int_0^{\infty}\frac{J_{\nu}(bx)x^{\nu+1}}{(x^2+a^2)^{\mu+1}}dx=\frac{a^{\nu-\mu}b^{\mu}}{2^{\mu}\Gamma(\mu+1)}K_{\nu-\mu}(ab),
\end{equation}
with $\nu=\frac{n-2}{2}, x=r, b= ||{\bf t}||, a=\sqrt{s}$,
we get
$$\int_0^{\infty}r^{\frac{n}{2}}J_{\frac{n-2}{2}}(r ||{\bf t}||)g(||{\bf t}||^2)dr=\frac{s^{\frac{3m+4n}{4}} ||{\bf t}||^{\frac{m+n-2}{2}}}
{2^{\frac{m+n-2}{2}}\Gamma(\frac{m+n}{2})}K_{-\frac{m}{2}}(\sqrt{s}||{\bf t}||),$$
where $K$ is the  modified Bessel function of the second kind.
It is easy to verify that
\begin{equation}
\int_0^{\infty}z^{\frac{n}{2}-1} g(z)dz=s^{\frac{n}{2}}\frac{\Gamma(\frac{n}{2})\Gamma(\frac{m}{2})}{\Gamma(\frac{n+n}{2})}.
\end{equation}
Substituting (3.9) and (3.10) into (3.8) and note that $K_{\nu}=K_{-\nu}$, we find that
\begin{equation}
\phi(u^2)=\frac{||{\bf t}||^{\frac{m}{2}}s^{\frac{m}{4}}}{2^{\frac{m}{2}-1}\Gamma(\frac{m}{2})}K_{\frac{m}{2}}(\sqrt{s}u), u\ge 0,
\end{equation}
which has been found by Joarder and Ali (1996), Joarder and Alam (1995) obtained the
characteristic function of the elliptic $t$-distribution using the conditional expectation technique, Song et al. (2014) got the same result by using the fact that the multivariate/generalized $t$ distributions can be expressed as a normal variance-mean mixture. For the special case of $n=1, m=s$, see Gaunt (2021).
In particular, when $m=1,s=1$, by using
$$K_{\frac12}(x)=\sqrt{\frac{\pi}{2x}}e^{-x},$$ one gets
the characteristic function of  multivariate Cauchy given by
\begin{equation}
\psi_{\bf X}({\bf t})=e^{i{\bf t}'{\boldsymbol \mu}}\exp(-\sqrt{{\bf t}'{\bf \Sigma}{\bf t}}),\; {\bf t}\in \Bbb{R}^n.
\end{equation}

\subsection{Pearson Type II }

 The $n$-dimensional random vector ${\bf X}$ is said to have a symmetric multivariate
Pearson type II distribution with parameters $m>-1$, ${\boldsymbol \mu}\in \Bbb{R}^n$, ${\bf \Sigma}: n\times n$ with ${\bf \Sigma}>0$  if its probability density function  is given by
$$f({\bf x})=\frac{\Gamma(\frac{n}{2}+m+1)}{\pi^{\frac{n}{2}}\Gamma(m+1)\sqrt{|{\bf \Sigma}|}}\left(1-({\bf x}- {\boldsymbol \mu})^{T}{\bf \Sigma}^{-1}({\bf x}-{\boldsymbol \mu})\right)^m,$$
where $0\le ({\bf x}- {\boldsymbol \mu})^{T}{\bf \Sigma}^{-1}({\bf x}-{\boldsymbol \mu})\le 1$.

This distribution was introduced by Kotz (1975) and will be denoted by $MPII_n({\boldsymbol \mu}, {\bf \Sigma})$.  A closed form of the characteristic function  of the multivariate Pearson type II distribution has been obtained by Joarder (1997).

By (3.2) and using the following integral (Gradshteyn and Ryzhik (2007, p.679))
$$\int_0^1 x^{\nu+1}(1-x^2)^{\mu}J_{\nu}(bx)dx=2^{\mu}\Gamma(\mu+1)b^{-(\mu+1)}J_{\nu+\mu+1}(b),$$
with $\mu=m$ and $\nu=\frac{n-2}{2}$, we get
the characteristic function  of  multivariate Pearson type II distribution:
\begin{equation}
\psi_{\bf X}({\bf t})=e^{i{\bf t}'{\boldsymbol \mu}}2^{\frac{n}{2}+m}\Gamma\left(\frac{n}{2}+m+1\right)(||{\bf \Sigma}^{\frac{1}{2}}{\bf t}||)^{-\frac{n}{2}-m}
J_{\frac{n}{2}+m}(||{\bf \Sigma}^{\frac 12}{\bf t}||),\; {\bf t}\in \Bbb{R}^n.
\end{equation}
 Using of the following relation between the Bessel
function of the first kind and the generalized hypergeometric function (Slater (1960), 1.8.5)
\begin{equation}
{}_0F_1\left(b+1,-\frac{x^2}{4}\right)=\left(\frac{x}{2}\right)^{-b}\Gamma(b+1)J_b(x),
\end{equation}
we get the following equivalent form of (4.10):
\begin{equation}
\psi_{\bf X}({\bf t})=e^{i{\bf t}'{\boldsymbol \mu}}{}_0F_1\left(\frac{n}{2}+m+1,-\frac{{\bf t}'{\bf \Sigma}{\bf t}}{4}\right), \; {\bf t}\in \Bbb{R}^n.
\end{equation}
The last expression has been obtained by Li (1994).
\begin{remark}  Result (3.15)  simplified the expression (2.2)-(2.4)  in Sutradhar (1986), and also
  revised   Joarder (1997, (2.1)) in which the following wrong identity was used:
\begin{equation}
{}_0F_1\left(b,-\frac{x^2}{4}\right)=\left(\frac{x}{2}\right)^{-b}\Gamma(b+1)J_b(x).
\end{equation}
\end{remark}

\subsection{Pearson Type VII }

Let ${\bf X}$ be an $n$-dimensional vector distributed according to a symmetric multivariate Pearson type VII
distribution with density  (cf. Fang et al. (1990)):
$$f({\bf x})=c_n\left(1+\frac{1}{m}({\bf x}- {\boldsymbol \mu})^{T}{\bf \Sigma}^{-1}({\bf x}-{\boldsymbol \mu})\right)^{-N},$$
where $N>\frac{n}{2}, m>0$ are parameters, $c_n$ is the normalizing constant given by
$$c_n=\frac{\Gamma(N)}{(m\pi)^{n/2}\Gamma(N-n/2)}.$$
This subclass includes a number of important distributions such as   the multivariate
$t$-distribution for $N=\frac{n+ m}{2}$ and the multivariate Cauchy distribution
for $m=1$ and   $N=\frac{n+ 1}{2}$.

Using  the following integral formula (see Gradshteyn and Ryzhik (2007, p.678))
$$
\int_0^{\infty}\frac{J_{\nu}(bx)x^{\nu+1}}{(x^2+a^2)^{\mu+1}}dx=\frac{a^{\nu-\mu}b^{\mu}}{2^{\mu}\Gamma(\mu+1)}K_{\nu-\mu}(ab),
$$
we get
\begin{equation}
\int_0^{\infty}r^{\frac{n}{2}}J_{\frac{n-2}{2}}(r ||{\bf t}||)\left(\frac{m+r^2}{m}\right)^{-N}dr
=\frac{m^{\frac{n}{4}+\frac{N}{2}}||{\bf t}||^{N-1}}{2^{N-1}\Gamma(N)}K_{\frac{n}{2}-N}(\sqrt{m}||{\bf t}||).
\end{equation}
It is easy to verify that
\begin{equation}
\int_0^{\infty}r^{\frac{n}{2}-1}\left(1+\frac{r}{m}\right)^{-N}dr=m^{\frac{n}{2}}\frac{\Gamma(\frac{n}{2})\Gamma(N-\frac{n}{2})}{\Gamma(N)}.
\end{equation}
Substituting (3.17) and (3.18) into (3.2) we obtain
\begin{equation}
\psi(||{\bf t}||^2)=e^{i{\bf t}'{\boldsymbol \mu}}\phi({\bf t}'{\bf \Sigma}{\bf t}),\; {\bf t}\in \Bbb{R}^n,\nonumber
\end{equation}
and
 \begin{eqnarray}
 \phi(u^2)&=&\frac{2^{\frac{n}{2}-N+1}}{\Gamma(N-\frac{n}{2})}m^{\frac{N}{2}-\frac{n}{4}}u^{N-\frac{n}{2}}K_{\frac{n}{2}-N}(\sqrt{m}u)\nonumber\\
 &=&\frac{2^{\frac{n}{2}-N+1}}{\Gamma(N-\frac{n}{2})}m^{\frac{N}{2}-\frac{n}{4}}u^{N-\frac{n}{2}}K_{N-\frac{n}{2}}(\sqrt{m}u).
 \end{eqnarray}
In particular, when $N=(n+m)/2$ and $m=s$, we recover the result (3.11).

The   characteristic  generator of ${\bf X}$ is obtained in   Fang et al. (1990, (3.30)):
$$ \phi(u^2)=\frac{2\Gamma(N-(n-1)/2)}{\sqrt{\pi}\Gamma(N-n/2)}\int_0^{\infty} \cos(\sqrt{m}ut)(1+t^2)^{-N+(n-1)/2}dt,$$
which is equivalent to (3.19). In fact,
by using the following integral (see Gradshteyn and Ryzhik (2007, p.442))
$$\int_0^{\infty}(\beta^2+x^2)^{\nu-1/2}\cos(ax)dx=\frac{1}{\sqrt\pi}\cos(\pi \nu)\left(\frac{2\beta}{a}\right)^{\nu}\Gamma(\nu+\frac12)K_{-\nu}(a\beta),$$
where $a>0, \Re(\beta)>0, \Re(\nu)<\frac12,$
we can rewritten above $\phi$ as
\begin{eqnarray*}
\phi(u^2)&=&\frac{2\Gamma(N-(n-1)/2)}{\sqrt{\pi}\Gamma(N-n/2)}\frac{1}{\sqrt\pi}\cos\left(\pi(\frac{n}{2}-N)\right)
\left(\frac{2}{\sqrt{m}u}\right)^{n/2-N}\\
&&\times\Gamma\left(\frac{n}{2}-N+\frac12\right)K_{N-n/2}(\sqrt{m}u)\\
&=&\frac{2^{\frac{n}{2}-N+1}}{\Gamma(N-\frac{n}{2})}\left(\frac{1}{\sqrt{m}u}\right)^{\frac{n}{2}-N}K_{N-\frac{n}{2}}(\sqrt{m}u),
 \end{eqnarray*}
 where we have used the following formula for gamma function:
  $$\Gamma(x)\Gamma(1-x)=\frac{\pi}{\sin\pi x}, x\in(0,1),$$
  as desired.

\subsection{Kotz type }

 Let ${\bf X}$ be distributed
according to a symmetric Kotz type distribution with density
\begin{equation}
f({\bf x})=c_n|{\bf \Sigma}|^{-\frac{1}{2}}[({\bf x}- {\boldsymbol \mu})^{T}{\bf \Sigma}^{-1}({\bf x-{\boldsymbol \mu} })]^{N-1}
\exp\left(-r[({\bf x}- {\boldsymbol \mu})^{T}{\bf \Sigma}^{-1}({\bf x-{\boldsymbol \mu} })]^s\right),
\end{equation}
for ${\bf x}\in \Bbb{R}^n$,  where ${\boldsymbol \mu}\in \Bbb{R}^n$, ${\bf \Sigma}$ is a positive  $n\times n$ matrix,
 $r,s >0, 2N+n>2$ are parameters and  $c_n$ is a normalizing constant:
 $$c_n=\frac{s\Gamma((\frac{n}{2})}{\pi^{\frac{n}{2}}\Gamma(2N+n-2)/2s)}r^{(2N+n-2)/2s}.$$
 The density generator takes the form
 $$g(l)=l^{N-1}e^{-rl^s}, \;l>0.$$ We denote it by ${\bf X}\sim KTD_n({\boldsymbol \mu}, {\bf \Sigma}, N, r,s).$
When $N = s = 1$ and $r = \frac12$, the distribution reduces to a multivariate
normal distribution and when $N = 1$ and $r =\frac12$ the distribution reduces
to a multivariate power exponential distribution.  A series form of the characteristic function  of the  symmetric Kotz type distribution has been obtained by  Iyengar and Tong (1989), Li (1994). Here we give an integral form in terms of the Bessel functions.
In doing so,  by (3.1) and (3.2) the characteristic function of ${\bf X}$  is given by
\begin{equation}
\psi_{\bf X}({\bf t})=e^{i{\bf t}'{\boldsymbol \mu}}\phi({\bf t}'{\bf \Sigma}{\bf t}),\; {\bf t}\in \Bbb{R}^n,
\end{equation}
 where
 \begin{eqnarray}
\phi(u^2)&=&c_n(2\pi)^{\frac{n}{2}} u^{-\frac{n-2}{2}}\int_0^{\infty}l^{\frac{n}{2}+2(N-1)}J_{\frac{n-2}{2}}(lu)e^{-rl^{2s}}dl\nonumber\\
&=&\frac{2^{\frac{n}{2}}s\Gamma((\frac{n}{2})r^{(2N+n-2)/2s}}{\Gamma(2N+n-2)/2s)} u^{-\frac{n-2}{2}}\nonumber\\
&&\times \int_0^{\infty}l^{\frac{n}{2}+2(N-1)}J_{\frac{n-2}{2}}(lu)e^{-rl^{2s}}dl, \; u\ge 0.
 \end{eqnarray}
In particular, when $s=1$, using the following integral (see Gradshteyn and Ryzhik (2007, p.706))
$$\int_0^{\infty} x^{\mu} e^{-\alpha x^2}J_{\nu}(\beta x)dx=\frac{\beta^{\nu}\Gamma(\frac{1}{2}(\nu+\mu+1))}{2^{\nu+1}\alpha^{\frac{1}{2}(\nu+\mu+1)}\Gamma(\nu+1)}
{}_1F_1\left(\frac{1}{2}(\nu+\mu+1);\nu+1;-\frac{\beta^2}{4\alpha}\right),$$
we get
\begin{eqnarray}
\phi(u^2)={}_1F_1\left(\frac{n}{2}+N-1;\frac{n}{2};-\frac{u^2}{4r}\right), \; u\ge 0,
 \end{eqnarray}
 which is the same form as of Li (1994);
when $s=\frac12$, using the integral (see Gradshteyn and Ryzhik (2007, p.702))
$$\int_0^{\infty}e^{-\alpha x}J_{\nu}(\beta x) x^{\nu+1}dx=\frac{2\alpha(2\beta)^{\nu}\Gamma(\nu+\frac{3}{2})}{\sqrt{\pi}(\alpha^2+\beta^2)^{\nu+\frac{3}{2}}},$$
we get
\begin{eqnarray}
\phi(u^2)=\frac{2^{n-1}\Gamma(\frac{n}{2})\Gamma(\frac{n+1}{2}) r^{2N+n-1}}
{\sqrt{\pi}\Gamma(2N+n-2) (r^2+u^2)^{\frac{n+1}{2}}}, \; u\ge 0.
 \end{eqnarray}
When $n=2$ and $N=1$, (3.24) reduces to the simple form (see Nadarajah (2003))
\begin{eqnarray}
\phi(u^2)=\frac{r^{3}}
{(r^2+u^2)^{\frac{3}{2}}}, \; u\ge 0.
 \end{eqnarray}

\subsection{Multivariate Bessel }

An $n\times 1$ random vector ${\bf X}$ is said to have a symmetric multivariate
Bessel distribution if the density has the form (see, e.g.  Fang et al. (1990))
\begin{equation*}
f({\bf x})=C_n|{\bf \Sigma}|^{-\frac{1}{2}}g(({\bf x}- {\boldsymbol \mu})^{T}{\bf \Sigma}^{-1}({\bf x-{\boldsymbol \mu} })), \;{\bf x}\in \Bbb{R}^n,
\end{equation*}
where
$$g(t)=\left(\frac{\sqrt{t}}{\beta}\right)^a K_a\left(\frac{\sqrt{t}}{\beta}\right),\; a>-\frac{n}{2}, \beta>0,$$
$$C_n^{-1}=2^{a+n-1}\pi^{\frac{n}{2}}\beta^{n+a}\Gamma\left(a+\frac{n}{2}\right).$$
It can be shown that
$$\int_0^{\infty}r^{\frac{n}{2}}J_{\frac{n-2}{2}}(ru)r^a K_a\left(\frac{r}{\beta}\right)dr
=2^{\frac{n}{2}+a-1}u^{\frac{n}{2}-1}\beta^{-a}\frac{\Gamma\left(\frac{n}{2}+a\right)}{(u^2+\beta^{-2})^{\frac{n}{2}+a}}, $$
from which and using (3.2) we get that
$$\phi(u^2)=\frac{1}{(1+u^2\beta^2)^{\frac{n}{2}+a}}, u\ge 0.$$
So the characteristic function of ${\bf X}$  is given by
\begin{equation*}
\psi_{\bf X}({\bf t})=\frac{e^{i{\bf t}'{\boldsymbol \mu}} }{(1+\beta^2{\bf t}'{\bf \Sigma}{\bf t} )^{\frac{n}{2}+a}},\; {\bf t}\in \Bbb{R}^n.
\end{equation*}
In particular, when $a=\frac12$ and $\beta=1$, we get the characteristic function of the Type-II multivariate Laplace distribution (see, Shi et al. (2022)):
\begin{equation*}
\psi_{\bf X}({\bf t})=\frac{e^{i{\bf t}'{\boldsymbol \mu}} }{(1+{\bf t}'{\bf \Sigma}{\bf t} )^{\frac{n+1}{2}}},\; {\bf t}\in \Bbb{R}^n.
\end{equation*}

\subsection{Multivariate Logistic}

An elliptical vector ${\bf X}$  belongs to the family of multivariate logistic distributions if its density
generator has the form
$$g(u)=\frac{e^{-u}}{(1+e^{-u})^2}.$$
A generalized elliptically symmetric logistic  distribution
including the multivariate   logistic distribution can be found in Yin et al. (2022).
By using (3.1) and (3.2) we find the characteristic function of ${\bf X}$ given by
\begin{equation}
\psi_{\bf X}({\bf t})=e^{i{\bf t}'{\boldsymbol \mu}}\phi({\bf t}'{\bf \Sigma}{\bf t}),\; {\bf t}\in \Bbb{R}^n\nonumber
\end{equation}
 where
\begin{equation}
\phi(u^2)=c_n(2\pi)^{\frac{n}{2}} u^{-\frac{n-2}{2}}\int_0^{\infty}r^{\frac{n}{2}}J_{\frac{n-2}{2}}(ru) \frac{e^{-r^2}}{(1+e^{-r^2})^2}dr, \; u\ge 0.\nonumber
\end{equation}
Here, $c_n$ is the normalizing constant given by
$$c_n=\frac{\Gamma(n/2)}{(2\pi)^{n/2}}\left(\int_0^{\infty} x^{n/2-1}\frac{e^{-x}}{(1+e^{-x})^2}dx\right)^{-1}.$$
The   characteristic  generator of ${\bf X}$ is given in Li (1994) by
$$\phi(x)=\sum_{k=0}^{\infty}\left(\frac{x}{4\pi}\right)^k\frac{c_n}{k!c_{n+2k}}.$$

\subsection{ Multivariate exponential power}
An elliptical vector ${\bf X}$  is said to have a multivariate
exponential power distribution if its density generator
has the form
$$g(u)=\exp(-ru^s),\; r>0,s>0.$$
When $r=s=1$,  this family of distributions clearly
reduces to the multivariate normal family. If $s=\frac12, r=\sqrt{2}$, we have the family of double
exponential or Laplace distributions.
By using (3.1) and (3.2) we find the characteristic function of ${\bf X}$ given by
\begin{equation}
\psi_{\bf X}({\bf t})=e^{i{\bf t}'{\boldsymbol \mu}}\phi({\bf t}'{\bf \Sigma}{\bf t}),\; {\bf t}\in \Bbb{R}^n\nonumber
\end{equation}
 where
\begin{equation}
\phi(u^2)=c_n(2\pi)^{\frac{n}{2}} u^{-\frac{n-2}{2}}\int_0^{\infty}r^{\frac{n}{2}}J_{\frac{n-2}{2}}(ru)\exp(-ru^{2s}) dr, \; u\ge 0.\nonumber
\end{equation}
Here, $c_n$ is the normalizing constant given by
\begin{eqnarray*}
c_n&=&\frac{\Gamma(n/2)}{(2\pi)^{n/2}}\left(\int_0^{\infty} x^{n/2-1} \exp(-rx^s) dx\right)^{-1}\\
&=& \frac{s\Gamma(n/2)}{(2\pi)^{n/2} \Gamma(n/2s)} r^{n/2s}.
\end{eqnarray*}

\section{Some  extensions}
In this section, we  consider  CFs of  the location-scale mixture of elliptical distributions and CFs of multivariate generalized skew-elliptical distribution, which can be as  the  generalizations of  elliptical distributions.

\subsection{CFs of location-scale mixture of elliptical distributions}
An $n$-dimensional random variable ${\bf X}$ is said to have a  location-scale mixture of elliptical distributions with the parameters ${\boldsymbol \mu}, {\boldsymbol\gamma}$ and  ${\bf \Sigma}$, if
\begin{equation}
{\bf X}={\boldsymbol \mu}+V{\boldsymbol \gamma}+\sqrt{V}{\bf \Sigma}^{\frac12}{\bf Z},
\end{equation}
where ${\bf Z} \sim ELL_n ({\bf 0},{\bf I}_n,g)$, $V$ is a nonnegative, scalar-valued random variable with the distribution $F$, ${\bf Z}$ and $V$ are independent,    ${\boldsymbol \mu}, {\boldsymbol \gamma}\in \Bbb{R}^n$,  ${\bf \Sigma}\in\Bbb{R}^{n\times n}$ with ${\bf \Sigma}>0$, and  ${\bf \Sigma}^{\frac12}$ is the square root of  ${\bf \Sigma}$. Here ${\bf 0}$ is an $n\times1$ vector of zeros, and ${\bf I}_n$ is $n\times n$ identity matrix.

Note that when ${\bf Z} \sim KTD_n({\bf 0},{\bf I_n}, N,\frac12, s)$ we have  the variance-mean mixture of the Kotz-type distribution  introduced by Arslan (2009); When ${\bf Z} \sim N_n(0,{\bf I_n})$ we get the multivariate normal variance-mean mixture distribution (see, e.g., McNeil et al.
(2005)).

\begin{theorem}
   The characteristic function of ${\bf X}$  defined by (4.1) has the form
\begin{equation}
\psi_{\bf X}({\bf t})=c_n (2\pi)^{\frac{n}{2}}e^{i{\bf t}'{\boldsymbol \mu}} \int_0^{\infty}e^{iv{\bf t}'{\boldsymbol \gamma}} \phi(v{\bf t}'{\bf \Sigma}{\bf t})dF(v) ,\; {\bf t}\in \Bbb{R}^n,
\end{equation}
where $c_n$ is defined by (2.3) and  $\phi$ is defined by (3.2).
\end{theorem}
{\bf Proof}\,  By using (3.1) and (3.2),
the characteristic function of  ${\bf X}$   can be written as
\begin{eqnarray}
\psi_{\bf X}({\bf t})&=& e^{i{\bf t}'{\boldsymbol \mu}} E\{E(e^{i{\bf t}'{(V{\boldsymbol \gamma}+\sqrt{V}{\bf \Sigma}^{\frac12}{\bf Z})}}|V)\}\nonumber\\
&=&e^{i{\bf t}'{\boldsymbol \mu}} E\{e^{iV{\bf t}'{\boldsymbol \gamma}} E(e^{i{\bf t}'{\sqrt{V}{\bf \Sigma}^{\frac12}{\bf Z}}}|V)\}\nonumber\\
&=&e^{i{\bf t}'{\boldsymbol \mu}} \int_0^{\infty}e^{iv{\bf t}'{\boldsymbol \gamma}} E(e^{i{\bf t}'{\sqrt{v}{\bf \Sigma}^{\frac12}{\bf Z}}})P(V\in dv) \nonumber\\
&=&e^{i{\bf t}'{\boldsymbol \mu}} \int_0^{\infty}e^{iv{\bf t}'{\boldsymbol \gamma}} \phi(v{\bf t}'{\bf \Sigma}{\bf t})P(V\in dv)  \nonumber\\
&=&e^{i{\bf t}'{\boldsymbol \mu}} \int_0^{\infty}e^{iv{\bf t}'{\boldsymbol \gamma}} \phi(v{\bf t}'{\bf \Sigma}{\bf t})dF(v). \nonumber
\end{eqnarray}
 This ends the proof of Theorem 4.1.

\subsection{CFs of  scale mixture of the uniform distributions}

A random vector ${\bf X}\sim S_n(g)$ is said to have a scale mixture of the uniform (SMU) representation, if it can be
represented as ${\bf X}=W{\bf V^{(n)}}$, where  ${\bf V^{(n)}}$ is  the random vector distributed uniformly inside a unit
sphere in $\Bbb{R}^n$, $W$ is a positive random variable and is distributed independently of ${\bf V^{(n)}}$.
The SMU class is more general than the normal scale-mixing  class and the SMU property is characterised by Theorem 2 in
 Fung and Seneta (2008). ${\bf X}=W{\bf V^{(n)}}$   if and only if ${\bf X}$ is star unimodal which is equivalent to $g'(x)\le 0$ for $x>0$. In this case, $W$ has a density given by
  \begin{equation}
  f_{W}(w)=-\frac{2\pi^{\frac{n}{2}}}{n\Gamma(\frac{n}{2})}w^n g'(w^2),\; w>0.
  \end{equation}
 Examples of  star unimodal elliptical distributions are  multivariate normal,   multivariate exponential power,  multivariate $t$,  multivariate Cauchy, Stable laws and Pearson type VII distributions, and so on.

 The characteristic function of ${\bf X}$  is given by
 \begin{eqnarray}
\psi_{\bf X}({\bf t})&=&E(E(e^{iW{\bf t}'{\bf V^{(n)}}}|W))\nonumber\\
&=&\int_0^{\infty}E(e^{iw{\bf t}'{\bf V^{(n)}}})P(W\in dw)\nonumber\\
&=&-(2\pi)^{\frac{n}{2}}||{\bf t}||^{-\frac{n}{2}}\int_0^{\infty}w^{\frac{n}{2}}J_{\frac{n}{2}}(w||{\bf t}||)g'(w^2)dw,\; {\bf t}\in \Bbb{R}^n,\nonumber
 \end{eqnarray}
where we have used (3.3) and (4.3).

 \subsection{CFs of multivariate generalized skew-elliptical distribution}

An $n$-dimensional random variable ${\bf Y}$ is said to have a multivariate generalized skew-elliptical (GSE) distribution, if its density has the following form (cf. Genton and Loperfido, 2005)
\begin{equation}
f({\bf x})=2|{\bf \Sigma}|^{-\frac{1}{2}}c_ng(({\bf x}- {\boldsymbol \mu})^{T}{\bf \Sigma}^{-1}({\bf x-{\boldsymbol \mu} }))
\pi({\bf \Sigma}^{-\frac{1}{2}}({\bf x-{\boldsymbol \mu} })), \;{\bf x}\in \Bbb{R}^n,
\end{equation}
where  ${\bf \Sigma}^{-\frac{1}{2}}$ is the inverse of ${\bf \Sigma}^{\frac{1}{2}}$,
$$|{\bf \Sigma}|^{-\frac{1}{2}}c_ng(({\bf x}- {\boldsymbol \mu})^{T}{\bf \Sigma}^{-1}({\bf x-{\boldsymbol \mu} }))$$
is a probability density function  of $n$-variate elliptical distribution ${\bf X} \sim ELL_n ({\boldsymbol \mu},{\bf \Sigma},g)$,
 $c_n$ is the normalizing constant given by (2.3), $\pi$ is the skewing function which satisfies $0\le \pi ({\bf y})\le 1$ and
$\pi (-{\bf y})+\pi ({\bf y})=1$ for ${\bf y}\in \Bbb{R}^n$.

 If  ${\bf \pi(y)}=F_{g^*}({\boldsymbol \alpha}'{\bf y})$, where $F_{g^*}(\cdot)$ is a univariate cdf of standard elliptical distribution with density generator $g^*$,  one gets the skew-elliptical  distribution introduced in Vernic (2006); If  ${\bf \pi(y)}=H({\boldsymbol \alpha}'{\Sigma}^{\frac{1}{2}}{\bf y})$, where $H(\cdot)$  is the cdf of a distribution symmetric around 0,
 one gets the skew-elliptical  distribution introduced in   Azzalini and  Dalla-Valle (1996).
The skew-elliptical  distributions   include the more familiar
skew-normal, skew-$t$ and skew-Cauchy distributions.

The characteristic function  of (4.4), in the special case of ${\boldsymbol \mu} =0$ and ${\bf \Sigma}={\bf I_n}$, has been derived by  Shushi (2016) as follows:
 \begin{equation}
 \psi_{\bf X}({\bf t})=2\phi({\bf t}'{\bf t})k_n({\bf t}), \;{\bf t}\in \Bbb{R}^n,
 \end{equation}
where $\phi$ is the characteristic generator of ${\bf X}$, $k_n$ is a function that satisfies
 $0\le k_n({\bf t})\le 1$ and
$k_n(-{\bf t})+k_n({\bf t})=1$ for ${\bf t}\in \Bbb{R}^n$.

In general case, $\psi_{\bf X}({\bf t})$ has the form (see,  Shushi (2018)):
 \begin{equation}
 \psi_{\bf X}({\bf t})=2e^{i{\bf t}'{\boldsymbol \mu}}\phi({\bf t}'{\bf \Sigma}{\bf t})a_n({\bf t}), \;{\bf t}\in \Bbb{R}^n,
 \end{equation}
 where
  $$a_n({\bf t})=\frac{c_n}{\phi({\bf t}'{\bf \Sigma}{\bf t})}\int_{\Bbb{R}^n}\cos({\bf t}'{\bf y})g({\bf y}'{\bf \Sigma}^{-1}{\bf y})\pi({\bf y})d{\bf y}.$$
Here, $a_n$  satisfies
 $0\le a_n({\bf t})\le 1$ and
$a_n(-{\bf t})+a_n({\bf t})=1$ for ${\bf t}\in \Bbb{R}^n$.

 Applying Theorem 3.1, the characteristic generator $\phi$ in (4.5) and (4.6) is given by
 \begin{equation}
\phi(u^2)=c_n(2\pi)^{\frac{n}{2}} u^{-\frac{n-2}{2}}\int_0^{\infty}r^{\frac{n}{2}}J_{\frac{n-2}{2}}(ru)g(r^2)dr, \; u\ge 0,\nonumber
\end{equation}
where $c_n$ is the normalizing constant given by (2.3).

\section{Conclusions}
In this paper, we  have  given a unified argument for  the characteristic functions of all elliptical  and  related
distributions in terms of  Bessel functions which avoid
the computation of the contour integration.  The approach was based on the stochastic representation of elliptical random variable
and the characteristic  function of
uniform distribution on the unit sphere surface in $\Bbb{R}^n$.  In particular, we present the  simple closed form of characteristic functions for  commonly used distributions such as  multivariate $t$,  Pearson Type II, Pearson Type VII, Kotz type  and  Bessel distributions.
The results can be easily generalized to the case
having location and scale parameters and to the skewed multivariate cases. We remark that the obtained results can be also expressed in terms of the  confluent hypergeometric  functions.

\vskip 0.2cm
\noindent{\bf Conflict of Interests.}\  The authors declare that there is no conflict of interests regarding the publication of this article.\\

\noindent{\bf\Large Acknowledgements} \  The authors thank the editor and three referees for helpful comments which improved an
earlier version of the work.
 This research was supported by the National Natural Science Foundation of China (No. 12071251).

\end{document}